\documentclass[aps,showkeys,showpacs,superscriptaddress,a4paper,twocolumn,10pt]{revtex4-2}%{article}

\usepackage{amssymb}
\usepackage{amsmath}
\usepackage[dvips]{graphicx}
\usepackage[T2A]{fontenc}
\usepackage{graphicx}
\usepackage[dvips]{color}
\usepackage[breaklinks=true,colorlinks=true,linkcolor=blue,urlcolor=blue,citecolor=blue]{hyperref}
\usepackage{subcaption}
\let\doi\relax
\usepackage{doi}

\begin{document}

\title[On the unconstrained Allen-Cahn system and its Cahn-Hilliard limit]{On the unconstrained Allen-Cahn system and its Cahn-Hilliard limit}

\author{P.~O.~Mchedlov-Petrosyan}
\author{L.~N.~Davydov}
\email[Corresponding author: ]{ldavydov@kipt.kharkov.ua}

\affiliation{National Science Center Kharkiv Institute of Physics and Technology, \\
1, Akademichna St., Kharkiv, 61108, Ukraine}

\author{O.~A.~Osmayev}
\affiliation{National Science Center Kharkiv Institute of Physics and Technology, \\
1, Akademichna St., Kharkiv, 61108, Ukraine}
\affiliation{Ukrainian State University of Railway Transport, 7, Feuerbach Sq., Kharkiv, 61050, Ukraine}

\begin{abstract}
In [\href{https://doi .org /10 .1007 /s11012 -021 -01338 -y}{F.P. Duda et al., Meccanica, 2021}] the new theory for constituent migration in solids with microstructure was developed. We consider the one-dimensional version of the basic system of this theory, with generally non-zero kinetic modulus and the particular diffusivity of the "fast diffusion" type. We obtain exact traveling wave solutions for this system, and analyze several unconstraint and constraint limiting cases. We show that in unconstraint cases of Allen-Cahn type, i.e. for two independent variables, the system demonstrates broader set of possible regimes; in the constraint limit only CH-type process remains.

\end{abstract}


\maketitle

\renewcommand{\theequation}{\arabic{section}.\arabic{equation}}
\section{Introduction}\label{s1}

In a recent article Duda, Sarmiento and Fried \cite{1} introduced an interesting generalizations of the well known Allen-Cahn \cite{2} and Cahn-Hilliard \cite{3,4} equations. The Allen-Cahn and Cahn-Hilliard models are, probably, the most popular models in the phenomenological theory of Phase Transitions and generally in Material Science. The phenomenological Allen-Cahn (from now on AC) theory presumes that for the evolution of the \textit{non-conserved} order parameter $w$ the driving force is the deviation of the chemical potential $\bar{\mu }$ from some equilibrium value, i.e.
\begin{equation} \label{1.1} \frac{\partial w}{\partial t'} =-R\bar{\mu } .  \end{equation}
Here $R$ is a kinetic coefficient. On the other hand the Cahn-Hilliard (from now on CH) theory presumes that the driving force for evolution of the \textit{conserved} order parameter $w$ is the gradient of the chemical potential, i.e.,
\begin{equation} \label{1.2} \frac{\partial w}{\partial t'} =\nabla \left(M\nabla \bar{\mu }\right) .  \end{equation}
For the inhomogeneous system the chemical potential is defined as the variational derivative of some thermodynamic potential functional $\Psi $
\begin{equation} \label{1.3} \bar{\mu }=\frac{\delta }{\delta w} \Psi  .  \end{equation}

The key-point of both theories is the idea (used already by Van der Waals in his theory of capillarity \cite{5}) that for the inhomogeneous system the thermodynamic potential functional should depend not only on the order parameter $w$, but on its derivatives as well, i.e., $\Psi =\Psi [\omega ,\nabla \omega ,\; \Delta \omega ,...]$.

In \cite{6,7} it was proposed an approach based ``on the separation of the basic balance laws (such as those for mass and force), which are general and hold for large classes of materials, from constitutive equations (such as those for elastic solids and viscous fluids), which delineate specific classes of material behavior'', see \cite{7}, p.~179. Then, introducing the balance equations for microforces and thermodynamically consistent constitutive equations, several generalizations of AC and CH equations were derived. This approach is further continued in \cite{1}: a new constrained theory for constituent migration in solids with microstructure was developed. The microstructure is described by a scalar phase field; the coincidence between the phase field and the constituent concentration is treated as an internal constraint that is maintained by appropriate reactions. ``This leads to an interpretation of the constrained theory underlying the Cahn--Hilliard family as a limiting variant of the unconstrained theory underlying the Allen--Cahn family'' \cite{1}, p.~1707.

Considering the parameter $\phi $, which characterizes the local ordering of atoms, and the mass fraction $c$ as independent variables, the authors derived corresponding general evolution equations. The free energy density is defined as
\begin{equation} \label{1.4} \psi =F\left(\phi \right)+\frac{\lambda }{2} \left|\nabla \phi \right|^{2} +\frac{\chi }{2} \left(c-\phi \right)^{2}  .  \end{equation}
Here $\chi >0$ is a constant coupling energy modulus. Then the chemical potential is
\begin{equation} \label{1.5} \bar{\mu }=\frac{\partial \psi }{\partial c} =\chi \left(c-\phi \right),   \end{equation}
and the particular form of the evolution equations for $c,\, \phi $ is
\begin{equation} \label{1.6} \frac{\partial c}{\partial t'} =\chi \nabla \left[\hat{K}\nabla \left(c-\phi \right)\right] ,  \end{equation}
\begin{equation} \label{1.7} \beta \left(\phi ,\frac{\partial \phi }{\partial t'} \, \right)\frac{\partial \phi }{\partial t'} =\lambda \Delta \phi -f\left(\phi \right)-\chi \left(\phi -c\right) .  \end{equation}
Here $f\left(\phi \right)=\frac{d}{d\phi } F\left(\phi \right)$. The left-hand side of \eqref{1.7} originates from the thermodynamically consistent general expression for the microforces density; the `kinetic modulus' $\beta \left(\phi ,\frac{\partial \phi }{\partial t'} \, \right)$ is a nonnegative function for each choice of its arguments. To further simplify the consideration, in \cite{1} it is presumed $\beta =0$ and $\hat{K}=k\hat{I}$, where $\hat{I}$ is a unit matrix.

In the present work we consider the one-dimensional version of the system \eqref{1.4}-\eqref{1.7}. We consider the mobility $K\left(c\right)$ corresponding to so-called `fast diffusion':
\begin{equation} \label{1.8} K=\frac{\bar{k}}{c} .  \end{equation}
The latter expression means that the mobility decreases with the increasing volume fraction of this component. It was first proposed in \cite{8}. This type of concentration dependence was observed in several systems, for instance in ion exchange kinetics \cite{9}, in thin liquid film dynamics \cite{10}, for the fluctuations dynamic in the self-organized systems \cite{11,12} and in Boltzman kinetic model \cite{13}. The diffusion equation with such dependence possesses some special features, which attracted much mathematical interest \cite{14,15,16,17,18,19}.

Also, we presume $\beta $ to be a positive constant and look for the consequences of such choice. Using \eqref{1.5} to eliminate $c$, we get
\begin{equation} \label{1.9} \frac{\partial \varphi }{\partial t'} +\frac{1}{\chi } \frac{\partial \bar{\mu }}{\partial t'} =\bar{k}\frac{\partial }{\partial x'} \left[\left(\varphi +\frac{\bar{\mu }}{\chi } \right)^{-1} \frac{\partial \bar{\mu }}{\partial x'} \right],  \end{equation}
\begin{equation} \label{1.10} \bar{\mu }=-\lambda \frac{\partial ^{2} \phi }{\partial x'^{2} } +f\left(\phi \right)+\beta \frac{\partial \phi }{\partial t'}  .  \end{equation}

If $\chi \to \infty $, in the terminology of \cite{1} the system becomes `constrained'; this formal limit for the $\beta =0;\, \, K=k$ case yields the standard CH system. This convergence was rigorously proven by Miranville \cite{20}. However, if $\beta >0$ and $K$ is given by \eqref{1.8}, the formal $\chi \to \infty $ limit reduces the system \eqref{1.9}-\eqref{1.10} to the \textit{viscous} Cahn-Hilliard equation \cite{21} with the viscosity $\beta $ and the non-constant mobility.

In the present paper solving the system \eqref{1.9}-\eqref{1.10} we obtain finally the exact traveling wave solution for the system \eqref{1.6}-\eqref{1.7}, i.e., $c$ and $\phi $. We show that solution for the `unconstrained' system exists only in a limited range of the coupling energy modulus $\chi $.

\setcounter{equation}{0}
\section{Traveling wave solution}\label{s2}

We consider the system \eqref{1.9}-\eqref{1.10} with generally asymmetric cubic polynomial $f\left(\phi \right)$ corresponding to the quartic polynomial $F\left(\phi \right)$:
\begin{equation} \label{2.1} \frac{\partial \varphi }{\partial t'} =\bar{k}\frac{\partial }{\partial x'} \left[\left(\varphi +\varepsilon \bar{\mu }\right)^{-1} \frac{\partial \bar{\mu }}{\partial x'} \right]-\varepsilon \frac{\partial \bar{\mu }}{\partial t'} ,  \end{equation}
\begin{equation} \label{2.2} \bar{\mu }=-\lambda \frac{\partial ^{2} \phi }{\partial x'^{2} } +\rho \phi ^{3} -\bar{\delta }\phi ^{2} -s\phi +\beta \frac{\partial \phi }{\partial t'}  .  \end{equation}
Here $\rho ,\, s$ are positive parameters and $\delta $ can be positive, negative, or zero; we also denoted $\frac{1}{\chi } =\varepsilon $. Introducing non-dimensional order parameter $u=\frac{\phi }{\phi _{0} }$  $\left(\phi _{0} =\frac{1}{\sqrt{\rho } } \right)$, non-dimensional coordinate $x=\frac{x'}{\sqrt{\lambda } } $, and non-dimensional time $t=\frac{t'}{T} $ $\left(T=\frac{\lambda }{\bar{k}\sqrt{\rho } } \right)$, we rewrite equations \eqref{2.1}, \eqref{2.2} as
\begin{equation} \label{2.3} \frac{\partial u}{\partial t} =\frac{\partial }{\partial x} \left[\frac{1}{\left(u+\varepsilon \mu \right)} \frac{\partial \mu }{\partial x} \right]-\varepsilon \frac{\partial \mu }{\partial t}  ,  \end{equation}
\begin{equation} \label{2.4} \mu =-\frac{\partial ^{2} u}{\partial x^{2} } +u^{3} -\delta u^{2} -su+\eta \frac{\partial u}{\partial t}  .  \end{equation}
Here we introduced notations
\begin{equation} \label{2.5} \mu =\, \frac{\bar{\mu }}{\phi _{0} } ,\, \, \, \eta =\frac{\beta }{T} ,\, \, \, \delta =\phi _{0} \bar{\delta }\, \,  .  \end{equation}

Looking for the traveling wave solution of \eqref{2.3}, \eqref{2.4} we introduce $z=x-vt$; then the system \eqref{2.3}-\eqref{2.4}  for $u\left(z\right),\, \mu \left(z\right)$ becomes
\begin{equation} \label{2.6} -v\frac{du}{dz} =\frac{d}{dz} \left[\frac{1}{\left(u+\varepsilon \mu \right)} \frac{d\mu }{dz} \right]+\varepsilon v\frac{d\mu }{dz} , \end{equation}
\begin{equation} \label{2.7} \mu =-\frac{d^{2} u}{dz^{2} } -v\eta \frac{du}{dz} +u^{3} -\delta u^{2} -su . \end{equation}

Integrating \eqref{2.6} once, we get
\begin{equation} \label{2.8} \frac{d\mu }{dz} =C\left(u+\varepsilon \mu \right)-v\left(u+\varepsilon \mu \right)^{2} . \end{equation}
Here $C$ is the integration constant. We look for a solution which approaches some constant values $u_{1} ,\, u_{2} $ if  $z\to \pm \infty $; for definiteness we presume $u_{1} <u_{2} $. For such solution the proper Ansatz is
\begin{equation} \label{2.9} \frac{du}{dz} =\kappa \left(u-u_{1} \right)\left(u-u_{2} \right)=\kappa \left(u^{2} -pu+q\right). \end{equation}
Here $p=u_{1} +u_{2} , \, \, q=u_{1} u_{2} $ and $\kappa $ is presently unknown constant. Using \eqref{2.9} expression \eqref{2.7} for the chemical potential $\mu $ is easily calculated
\begin{eqnarray} \label{2.10} & & {\mu =\left(1-2\kappa ^{2} \right)u^{3} +\left(3\kappa ^{2} p-v\eta \kappa -\delta \right)u^{2} }  \\ & & {+\left[-\kappa ^{2} \left(2q+p^{2} \right)+v\eta \kappa p-s\right]u+\left(\kappa ^{2} pq-v\eta \kappa q\right)}  . \nonumber \end{eqnarray}

Comparing powers of $u$ in the latter expression with the powers in \eqref{2.8} makes evident that the latter expression should be linear in $u$, i.e., $\mu =\gamma u+\zeta $. Here $\gamma ,\, \zeta $ are presently unknown constants. Then we can `split' equation \eqref{2.8} into two equalities
\begin{eqnarray} \label{2.11}  & &{\left(1-2\kappa ^{2} \right)u^{3} +\left(3\kappa ^{2} p-v\eta \kappa -\delta \right)u^{2} }  \\ & &{+\left[-\kappa ^{2} \left(2q+p^{2} \right)+v\eta \kappa p-s\right]u+\left(\kappa ^{2} pq-v\eta \kappa q\right)}\nonumber  \\ & &{=\gamma u+\zeta } , \nonumber  \end{eqnarray}
\begin{equation} \label{2.12} \gamma \frac{du}{dz} =C\left[u+\varepsilon \left(\gamma u+\zeta \right)\right]-v\left[u+\varepsilon \left(\gamma u+\zeta \right)\right]^{2}  .  \end{equation}
After some rearrangements Eq.~\eqref{2.12} becomes
\begin{eqnarray} \label{2.13} \gamma \frac{du}{dz} &=&-v\left(1+\varepsilon \gamma \right)^{2} u^{2} +\left[\left(1+\varepsilon \gamma \right)C-2\varepsilon \zeta v\left(1+\varepsilon \gamma \right)\right]u \nonumber \\ &+&\left[\varepsilon C\zeta -v\left(\varepsilon \zeta \right)^{2} \right]. \end{eqnarray}

Naturally, the latter equation should coincide with \eqref{2.9}; this yields three constraints on the parameters:
\begin{equation} \label{2.14} -v\left(1+\varepsilon \gamma \right)^{2} =\kappa \gamma  ,  \end{equation}
\begin{equation} \label{2.15} \left(1+\varepsilon \gamma \right)\left(C-2\varepsilon \zeta v\right)=-\kappa p\gamma  ,  \end{equation}
\begin{equation} \label{2.16} \varepsilon C\zeta -v\left(\varepsilon \zeta \right)^{2} =\kappa q\gamma  .  \end{equation}

Equality \eqref{2.11}  should be satisfied identically; this yields four constraints
\begin{equation} \label{2.17} 2\kappa ^{2} =1 ,  \end{equation}
\begin{equation} \label{2.18} 3\kappa ^{2} p-v\eta \kappa -\delta =0 ,  \end{equation}
\begin{equation} \label{2.19} -\kappa ^{2} \left(2q+p^{2} \right)+v\eta \kappa p-s-\gamma =0 ,  \end{equation}
\begin{equation} \label{2.20} \kappa ^{2} pq-v\eta \kappa q-\zeta =0 .  \end{equation}

If the constraints \eqref{2.14}-\eqref{2.20} are satisfied, the solution of \eqref{2.9} is simultaneously solution of \eqref{2.6}-\eqref{2.7}. Integrating and selecting $z=0$ to be the position of the $\max \left|\frac{du}{dz} \right|$ we get
\begin{equation} \label{2.21} u=\frac{u_{2} +u_{1} }{2} -\frac{u_{2} -u_{1} }{2} \tanh \left(\frac{1}{2} \kappa \left(u_{2} -u_{1} \right)\left(x-vt\right)\right). \end{equation}

\setcounter{equation}{0}
\section{Parametric dependence of solution for non-zero kinetic modulus}\label{s3}

The constraints \eqref{2.14}-\eqref{2.20} form a nonlinear algebraic system for seven unknowns $\kappa ,\, \, v,\, \, \gamma ,\, \, \zeta ,\, \, C,$ and $u_{1} ,\, u_{2} $ (or $p,\, q$). Solution of this system yields the dependence of the parameters of the solution \eqref{2.21} on the system parameters $\varepsilon ,\, \delta ,\, s$ and $\eta $. Eliminating $\kappa ^{2} $ (but leaving $\kappa $, it may be positive or negative) we transform \eqref{2.18}-\eqref{2.20} into
\begin{equation} \label{3.1} \frac{3}{2} p-v\eta \kappa -\delta =0 ,  \end{equation}
\begin{equation} \label{3.2} -q-\frac{1}{2} p^{2} +v\eta \kappa p-s-\gamma =0 ,  \end{equation}
\begin{equation} \label{3.3} \frac{1}{2} pq-v\eta \kappa q-\zeta =0 .  \end{equation}
Eliminating $v$ from \eqref{3.1}-\eqref{3.3}, we get
\begin{equation} \label{3.4} q-p^{2} +\delta p+s+\gamma =0 ,  \end{equation}
\begin{equation} \label{3.5} \left(p-\delta \right)q+\zeta =0 .  \end{equation}
Elimination of $C$ from \eqref{2.15}-\eqref{2.16}  yields
\begin{equation} \label{3.6} \, \left(1+\varepsilon \gamma \right)\left[\kappa q\gamma \, -\left(\varepsilon \zeta \right)^{2} v\right]=-\kappa p\gamma \varepsilon \zeta  .  \end{equation}
Elimination of $v$ from the latter equation and \eqref{2.14} yields
\begin{equation} \label{3.7} \left[q\left(1+\varepsilon \gamma \right)^{2} \, +\left(\varepsilon \zeta \right)^{2} +\varepsilon \zeta \left(1+\varepsilon \gamma \right)p\right]\gamma =0 .  \end{equation}
Eliminating $v$ from \eqref{2.14} and \eqref{3.1} we get
\begin{equation} \label{3.8} \, \left(\, 2\delta -3p\right)\left(1+\varepsilon \gamma \right)^{2} =\gamma \eta  .  \end{equation}

Equations \eqref{3.4}, \eqref{3.5} and \eqref{3.7}, \eqref{3.8} constitute a nonlinear algebraic system for four unknowns $p,\, q$ and $\gamma ,\, \zeta $. Evidently, \eqref{3.7} is satisfied if $\gamma =0$; however, it follows from \eqref{2.14} that $v=0$, i.e., only static solution is possible. Denoting the corresponding values of the parameters by the lower index `$s$', we have from \eqref{3.4}, \eqref{3.5} and \eqref{3.8}:
\begin{equation} \label{3.9} q_{s} -p_{s}^{2} +\delta p_{s} +s=0 ,  \end{equation}
\begin{equation} \label{3.10} \left(p_{s} -\delta \right)q_{s} +\zeta _{s} =0 ,  \end{equation}
\begin{equation} \label{3.11} \, p_{s} =\frac{2}{3} \delta  .  \end{equation}

Using \eqref{3.9}-\eqref{3.11} the values $u_{1s} ,\, u_{2s} $ of the order parameter and the values of the concentration $c_{1s} ,\, c_{2s} $ at $\pm \infty $ can be calculated, see Appendix 1. The $u_{1s} $ is always negative, however, the order parameter should not be necessary positive. On the other hand, only positive concentrations make sense. It appears that $c_{1s} >0$ only for negative $\delta $ and
\begin{equation} \label{3.12} \varepsilon \left|\delta \right|\left(\frac{2}{9} \delta ^{2} +s\right)\ge \sqrt{3\delta ^{2} +9s} +\left|\delta \right| .  \end{equation}

The latter inequality is violated if $\varepsilon \to 0$, i.e., the smaller boundary concentration becomes negative. Thus, the physical static solution exists only for the unconstraint system with the asymmetric potential, $\delta <0$.
If inequality \eqref{3.12} is satisfied, there are static solutions for $u$ and $c$:
\begin{equation} \label{3.13} u_{s} =\frac{u_{2s} +u_{1s} }{2} -\frac{u_{2s} -u_{1s} }{2} \tanh \left(\frac{1}{2} \kappa \left(u_{2s} -u_{1s} \right)x\right) ,  \end{equation}
\begin{equation} \label{3.14} c_{s} =\frac{1}{\sqrt{\rho } } \left[u_{s} +\varepsilon \frac{1}{3} \left|\delta \right|\left(\frac{2}{9} \delta ^{2} +s\right)\right] .  \end{equation}
For $\gamma \ne 0$ the constraint \eqref{3.7} becomes
\begin{equation} \label{3.15} q\left(1+\varepsilon \gamma \right)^{2} \, +\left(\varepsilon \zeta \right)^{2} +\varepsilon \zeta \left(1+\varepsilon \gamma \right)p=0 .  \end{equation}

The nonlinear algebraic system \eqref{3.4}, \eqref{3.5}, \eqref{3.8} and \eqref{3.15} is complicated; we have found only two exactly solvable limiting cases. We start the analysis of this system with the constraint case $\varepsilon =0$  and denote the unknowns by the zero lower index. Equation \eqref{3.15} is reduced to $q_{0} =0$, i.e., $u_{10} =0, \, \, p_{0} =u_{20} $ and equation \eqref{3.5} is reduced to $\zeta _{0} =0$. Equations \eqref{3.4} and \eqref{3.8} become
\begin{equation} \label{3.16} p_{0}^{2} -\delta p_{0} -s-\gamma _{0} =0 ,  \end{equation}
\begin{equation} \label{3.17} \, \, 2\delta -3p_{0} =\gamma _{0} \eta  .  \end{equation}
Here we are interested in the non-zero kinetic modulus $\beta $, i.e., $\eta \ne 0$. At the end of this Section~\ref{s3} we check the consequence of  $\eta =0$ in \eqref{3.17}.  The generally unconstraint case with zero kinetic modulus, i.e. $\varepsilon \ne 0;\, \, \eta =0$ is considered in the next Section~\ref{s4}. So now we presume $\eta \ne 0$. Then elimination of $\gamma _{0} $ yields the quadratic equation for $p_{0} $:
\begin{equation} \label{3.18} p_{0}^{2} -\left(\delta -\frac{3}{\eta } \right)p_{0} -\left(s+2\frac{\delta }{\eta } \right)=0 .  \end{equation}
The roots of the latter equation are
\begin{equation} \label{3.19} p_{0} =\frac{1}{2} \left(\delta -\frac{3}{\eta } \right)\pm \sqrt{\frac{1}{4} \left(\delta -\frac{3}{\eta } \right)^{2} +s\, +2\frac{\delta }{\eta } }  .  \end{equation}

For $\varepsilon =0$ we have $c_{20} =\varphi _{0} u_{20} $. Then to avoid the negativity of mobility \eqref{1.8} the value of the order parameter $u_{20} =p_{0} $ should be positive. If $s\eta +2\delta >0$ there is a single positive root, i.e., only `$+$' sign in \eqref{3.19} should be taken. So there is no solution for a strong asymmetry of the potential, $-2\delta >s\eta $. Correspondingly, $\gamma _{0} $ is

\begin{equation} \label{3.20} \gamma _{0} =\frac{9}{2\eta ^{2} } +\frac{1}{\eta } \left[\frac{\delta }{2} -3\sqrt{\frac{1}{4} \left(\delta -\frac{3}{\eta } \right)^{2} +s\, +2\frac{\delta }{\eta } } \right] .  \end{equation}

From the constraint \eqref{2.14} for $\varepsilon =0$ the velocity is
\begin{equation} \label{3.21} v_{0} =-\kappa \gamma _{0}.  \end{equation}
For the symmetric potential we have $\delta =0$, and \eqref{3.20} is simplified to
\begin{equation} \label{3.22} \gamma _{0} =\frac{9}{2\eta ^{2} } -\frac{3}{2\eta } \sqrt{\left(\frac{3}{\eta } \right)^{2} +4s\, }  .  \end{equation}
So in the constraint case and non-zero kinetic modulus the velocity of the transition front is non-zero even for the symmetric potential $\delta =0$ ($s$ is per definition positive).
Then from \eqref{3.19} it follows
\begin{equation} \label{3.23} u_{20} =p_{0} =-\frac{3}{2\eta } +\frac{1}{2} \sqrt{\left(\frac{3}{\eta } \right)^{2} +4s\, }.  \end{equation}
If the $\eta $ is small, $\, s\eta ^{2} \ll \frac{9}{4} $; then $\gamma _{0} \simeq -s$;  $v_{0} \simeq \kappa s$ and $\, u_{20} \simeq \frac{1}{3} \eta s$ ,
If the $\eta $  is large,  i.e., $\, s\eta ^{2} \gg \frac{9}{4} $, we have
\begin{equation}v_{0} =-\kappa \gamma _{0} \simeq \kappa \left(\frac{3\sqrt{s} }{\eta } -\frac{9}{2\eta ^{2} } +\frac{27}{8\eta ^{3} \sqrt{s} } \right), \label{3.24)}\end{equation}
\begin{equation}u_{20} \simeq \sqrt{s} -\frac{3}{2\eta } +\frac{9}{8\eta ^{2} \sqrt{s} } . \label{3.25)}\end{equation}
So, for the large viscosity $\eta $ (i.e., large kinetic modulus $\beta $) the velocity is decreasing with $\eta $, and  the concentration $u_{20} $ at $-\infty $ is close to the stationary state $\sqrt{s} $ of the symmetric potential.

For completeness we consider here the limit of zero kinetic modulus $\beta $, i.e., zero viscosity $\eta $. If $\eta =0$, equation \eqref{3.17} yields $p_{0} =\, u_{20} =\, \frac{2}{3} \delta $. We remind that $q_{0} =0$ and $u_{10} =0$; i.e., for the symmetric potential $\delta =0$ the wave amplitude is zero; for the positivity of concentrations it should be $\delta >0$. Equation \eqref{3.16} determines $\gamma _{0} $:
\begin{equation}\gamma _{0} =-\, \frac{2}{9} \delta ^{2} -s . \label{3.26}\end{equation}
Then the velocity of the transition front, see \eqref{3.21}, is
\begin{equation}v_{0} =\kappa \left(\frac{2}{9} \delta ^{2} +s\right) . \label{3.27}\end{equation}

\setcounter{equation}{0}
\section{Unconstraint case for zero kinetic modulus}\label{s4}

Another exactly solvable case of the system \eqref{3.4}, \eqref{3.5} and \eqref{3.8}, \eqref{3.15} is the unconstraint case $\varepsilon \ne 0$ with zero viscosity $\eta =0$. For zero viscosity Eq. \eqref{3.8} becomes
\begin{equation} \label{4.1} \, \left(\, 2\delta -3p\right)\left(1+\varepsilon \gamma \right)^{2} =0 . \end{equation}
We remind that in deriving \eqref{3.15} we presumed $\gamma \ne 0$; if $\, \left(1+\varepsilon \gamma \right)=0$ then Eq.~\eqref{2.14} could not be satisfied for the finite velocity $v$. Therefore, $\, \left(1+\varepsilon \gamma \right)\ne 0$ and the only possibility to satisfy \eqref{4.1} is
\begin{equation} \label{4.2} p=\frac{2}{3} \delta  .  \end{equation}
Then from \eqref{3.5}
\begin{equation} \label{4.3} q\delta =3\zeta  .  \end{equation}
If $\delta =0$, then $p=0, \, \, \zeta =0$; together with $\, \left(1+\varepsilon \gamma \right)\ne 0$ it follows from \eqref{3.15} $q=0$, i.e. zero amplitude wave. So it is necessary to have $\delta \ne 0$, i.e. asymmetric potential. For $\delta \ne 0$ \eqref{4.3} could be rewritten as
\begin{equation} \label{4.4} q=3\frac{\zeta }{\delta }  .  \end{equation}

Substituting \eqref{4.4} and \eqref{4.2} for the $q$ and $p$ into \eqref{3.15}  and introducing $\omega =1+\varepsilon \gamma $ for convenience, we get
\begin{equation} \label{4.5} \zeta \left(\omega ^{2} \, +\frac{\varepsilon ^{2} }{3} \delta \zeta +\varepsilon \omega \frac{2}{9} \delta ^{2} \right)=0 .  \end{equation}
If $\zeta =0$, then $q=0$, and  we recover the solution given at the end of the previous Section~\ref{s3}, i.e., $\varepsilon =0,\, \eta =0$ case. In the present Section~\ref{s4} we keep $\varepsilon \ne 0$; presuming $\zeta \ne 0$ we get
\begin{equation} \label{4.6} \omega ^{2} \, +\frac{\varepsilon }{9} \delta \left(3\varepsilon \zeta +2\delta \omega \right)=0 .  \end{equation}
Substitution of \eqref{4.2} and \eqref{4.4} into \eqref{3.4} yields
\begin{equation} \label{4.7} \, 3\varepsilon \zeta =-\delta \omega -\varepsilon \frac{2}{9} \delta ^{3} -\varepsilon s\delta +\delta  .  \end{equation}
Elimination of $\zeta $ from \eqref{4.6}-\eqref{4.7} yields the quadratic equation for $\omega $
\begin{equation} \label{4.8} \omega ^{2} \, +\frac{\varepsilon }{9} \delta ^{2} \omega -\frac{\varepsilon }{9} \delta ^{2} \left(\varepsilon \frac{2}{9} \delta ^{2} +\varepsilon s-1\right)=0 .  \end{equation}
The roots of the latter equation are
\begin{equation} \label{4.9} \omega _{1,2} =\frac{\varepsilon }{18} \left[-\delta ^{2} \pm 3\delta \sqrt{\left(\delta ^{2} +4s\right)-\frac{4}{\varepsilon } } \right] .  \end{equation}
To have the real value of $\omega $ the parameter $\varepsilon $ should be limited from below, $\varepsilon \ge \varepsilon _{l} $, where
\begin{equation} \label{4.10} \varepsilon _{l} =\frac{4}{\delta ^{2} +4s}  .  \end{equation}

So this solution exists only for the definitely unconstraint case, This is in agreement with the assumption $\varepsilon \ne 0$ which was essential in the derivation of equation \eqref{4.6}. Parameter $\omega $ given by \eqref{4.9} determines $\gamma $ and together with \eqref{2.14} the velocity of the transition fronts both of $u$ and $c$. Due to Eq.~\eqref{4.7} parameter $\zeta $ is also determined; consequently with \eqref{4.2} and \eqref{4.4} this yields $p$ and $q$ (or $u_{1} ,\, u_{2} $), i.e., the full parametric dependence of the solution \eqref{2.21}. However, the resulting exact expressions are too bulk and poorly informative. Considering the limiting case we presume $\varepsilon \gg \varepsilon _{l} $, i.e., relatively small $\chi $, and neglect the $\frac{4}{\varepsilon } $ term inside the radical in \eqref{4.9}. This situation may be called `strongly unconstraint'. For the further it is convenient to introduce the non-dimensional stable stationary states of the thermodynamic potential (see, (2.4)) $a_{1} ,\, a_{2} $  ($a_{1} <\, a_{2} $); then
\begin{equation} \label{4.11} \delta =a_{1} +a_{2} ,\, \, \, \, s=-a_{1} a_{2} . \end{equation}
The parameter $s$ is positive per definition, so $a_{1} <0$. Then the expression inside the radical in \eqref{4.9} could be rewritten as
\begin{equation} \label{4.12)} \delta ^{2} +4s=\left(a_{2} -\left|a_{1} \right|\right)^{2} +4a_{2} \left|a_{1} \right|=\left(a_{2} +\left|a_{1} \right|\right)^{2} . \end{equation}
Then if $\varepsilon \gg \varepsilon _{l} =\frac{4}{\left(\left|a_{1} \right|+a_{2} \right)^{2} } $ the roots \eqref{4.9} become
\begin{equation} \label{4.13} \omega _{1,2} =\frac{\varepsilon }{18} \delta \left[-\delta \pm 3\sqrt{\left(\delta ^{2} +4s\right)} \right] \end{equation}
or
\begin{equation} \label{4.14} \omega _{1} =\frac{\varepsilon }{9} \left(a_{2} -\left|a_{1} \right|\right)\left[a_{2} +2\left|a_{1} \right|\right] ,\end{equation}
\begin{equation} \label{4.15} \omega _{2} =-\frac{\varepsilon }{9} \left(a_{2} -\left|a_{1} \right|\right)\left[2a_{2} +\left|a_{1} \right|\right] .\end{equation}
Substitution of \eqref{4.14}, \eqref{4.15} into \eqref{4.7} yields two possible values of $\zeta $ and $q$, see \eqref{4.4}:
\begin{equation} \label{4.16} q_{1} =\, 3\frac{\zeta _{1} }{\delta } =-\frac{1}{3} a_{2} \left(a_{2} +2\left|a_{1} \right|\right), \end{equation}
\begin{equation} \label{4.17} q_{2} =3\frac{\zeta _{2} }{\delta } =-\frac{1}{3} \left|a_{1} \right|\left(2a_{2} +\left|a_{1} \right|\right) . \end{equation}

In \eqref{4.16}, \eqref{4.17} we have again neglected the term $\frac{1}{\varepsilon } $. Then there are two possible pairs $u_{1} ,\, u_{2} $, i.e. pairs of the roots of the equation
\begin{equation} \label{4.18} u^{2} -pu+q_{1,2} =0 ,  \end{equation}
corresponding to two different values of $q_{1,2} $; here the single value $p=\frac{2}{3} \delta $ is given by \eqref{4.2}. For $q_{1} $ given by \eqref{4.16} we denote the roots of \eqref{4.18} as $u_{11} ,\, u_{12} $; they are
\begin{equation} \label{4.19} u_{11} =-\frac{1}{3} \left(2\left|a_{1} \right|+a_{2} \right),\, \, \, u_{12} =a_{2} . \end{equation}
For the $q_{2} $ given by \eqref{4.17} the roots are
\begin{equation} \label{4.20} u_{21} =-\left|a_{1} \right|;\, \, \, u_{22} =\frac{1}{3} \left(2a_{2} +\left|a_{1} \right|\right). \end{equation}

To avoid the negativity of the concentration and the negativity of the diffusion coefficient \eqref{1.8}, the smaller concentration values $c_{11} ,\, \, c_{21} $ should be nonnegative. This is checked in the Appendix 2. The physical solutions exist only for $\delta <0$; for the both solutions the smallest values of concentrations are zero, $c_{11} =c_{21} =0$. For the both solutions the largest values of the concentrations are equal again, see, \eqref{7.5}, \eqref{7.6}
\begin{equation} \label{4.21} c_{12} =c_{22} =\frac{2\varepsilon }{27\sqrt{\rho } } \left|\delta \right|\left(a_{2} +2\left|a_{1} \right|\right)\left[\left|a_{1} \right|+2a_{2} \right]. \end{equation}

The limiting values are equal for the both solutions. However, while for the solution with $\omega _{2} $ the concentration and the order parameter increase/decrease simultaneously, for the solution with $\omega _{1} $ if the order parameter $u$ increases with $z$, the concentration $c$ decreases, and vice versa (see Appendix 2).

Neglecting again the $\frac{1}{\varepsilon } $ term, we get expression for the velocity of the transition fronts for the cases with $\omega _{2} $ and $\omega _{1} $
\begin{equation} \label{4.22} v=-\frac{9\kappa }{\varepsilon ^{2} \left|\delta \right|\left(2a_{2} +\left|a_{1} \right|\right)} , \end{equation}
\begin{equation} \label{4.23} v=\frac{9\kappa }{\varepsilon ^{2} \left|\delta \right|\left(a_{2} +2\left|a_{1} \right|\right)} . \end{equation}
For negative $\delta $ i.e., $\left|a_{1} \right|>a_{2} $, we have the inequality for denominators, $a_{2} +2\left|a_{1} \right|>2a_{2} +\left|a_{1} \right|$. Therefore, the absolute value of the velocity is larger for the case of \eqref{4.22}. For positive $\kappa $ the anti-kink for the order parameter moves in the negative direction in the case of \eqref{4.22}, and in the positive direction in the case of \eqref{4.23}. However, the concentration front is always moving diminishing the concentration.

\setcounter{equation}{0}
\section{Discussion}\label{s5}

There are several models with two order parameters, where one order parameter is concentration. In such models the thermodynamic potential depends both on the concentration and the order parameter; the chemical potentials for each parameter are defined as corresponding variational derivatives. Probably the first such model was introduced by Eguchi et al. \cite{22}.The evolution of the conserved order parameter (concentration) is governed by the CH equation, the evolution of the non-conserved order parameter (local ordering of atoms) by the Allen-Cahn equation. According to the general classification of Hohenberg and Halperin \cite{23} this is a type C model. In \cite{24} the connection between the discrete (lattice) approach and the continuum approach to the combination of conserved and non-conserved dynamics was studied in detail. After this work the name `Cahn-Hilliard/Allen-Cahn system' for these systems became commonly used. The two order parameters are generally independent variables; imposing a particular algebraic link between them results in a rather special solution \cite{25}.

On the other hand, in \cite{1} while the evolution of the order parameter (phase field) is governed by the Allen-Cahn-type equation \eqref{2.7}, the evolution of the concentration is governed by the diffusion equation with $-\nabla \mu $ as a usual driving force, see \eqref{1.6}. With the free energy density \eqref{1.4}, it was found on the basis of this approach that the CH equation can be interpreted as the limiting variant of this system for zero `kinetic modulus' $\beta $. This result was supported by numerical simulation.

In the terminology of \cite{1} the system \eqref{1.9}-\eqref{1.10} becomes constrained if $\chi \to \infty $. In this case the formal limit for $\beta =0 $, $\hat{K} = k\hat{I} $ results in the standard CH system. We have considered the case of generally non-zero kinetic modulus, $\beta >0$. In this case the formal constraint limit $\chi \to \infty $ results in the viscous CH equation. Different from the convective-viscous CH equation \cite{26}, there is no exact solution for the viscous CH equation with the constant mobility. Therefore, we have selected particular form of mobility \eqref{1.8} which allows the exact traveling wave solutions for the viscous CH equation. The physical circumstances in which the mobility expressions of this type are appropriate are discussed in the cited articles. They are rather different for different systems and depend on the particular diffusion mechanisms; the discussion of these mechanisms goes far beyond the scope of the present paper.

We have found exact traveling wave solutions for the unconstraint system \eqref{1.9}-\eqref{1.10} and its constrained limit. The method of solution is sometimes called `Riccati equation method'. The Ansatz \eqref{2.9} we used is a guess: the simplest possible expression for the derivative of solution with proper asymptotic behavior at $\pm \infty $. We use this Ansatz to calculate the chemical potential \eqref{2.7}, obtaining expression \eqref{2.10}. Comparing the powers of $u$ in the latter expression with the powers in \eqref{2.8} makes evident that the latter expression should be linear in $u$. This allows to `split' equation \eqref{2.8} into first order differential equality \eqref{2.12} and a polynomial in $u$ equality \eqref{2.11}. The differential equality evidently should coincide with the Riccati equation \eqref{2.9}; the polynomial equality should be fulfilled for arbitrary $u$. These requirements lead to the system of algebraic constraints imposed on the parameters of the system $\varepsilon,\, \delta,\, s, \, \eta $, velocity of the transition front $v$, parameters of the Ansatz \eqref{2.9}, i.e. $u_{1},\,\, u_{2}$ (or $p,\,\, q$), $\kappa $ and the auxiliary parameters $\gamma,\, \zeta $. If these constraints are fulfilled, the solution \eqref{2.21} of the Riccati equation \eqref{2.9} is a special solution for the system \eqref{2.6}-\eqref{2.7}. Then in Sections \ref{s3} and \ref{s4} the solution parameters are analyzed as functions of the system parameters for three exactly solvable cases of the algebraic system \eqref{2.14}-\eqref{2.20}.

The first solvable case corresponds to a static solution with zero velocity $v$ of the transition front. Formally there are always real limiting values for the order parameter and consequently the solution \eqref{2.21} is admissible. However, the lower limiting value for the concentration is positive only for asymmetric potential $\delta <0$ and large enough $\varepsilon $ (small coupling energy modulus $\chi $), see \eqref{3.12}. Therefore, the physically admissible solution exists in the strictly unconstraint case only.

The second case is the constrained limit $\varepsilon \to 0$ ($\chi \to \infty $) when the nondimensional concentration coincides with the order parameter $u$. In this limit the smaller limiting value of the order parameter is zero. To avoid the negativity of mobility the larger limiting value of the order parameter $u_{20} $ should be strictly positive. The positivity of $u_{20} $ is guaranteed if $2\left(a_{2} -\left|a_{1} \right|\right)+\left|a_{1} \right|a_{2} \eta >0$; the violation of this inequality means that there is no admissible solution. If the kinetic modulus $\beta =\eta T$ is non-zero, the velocity of the phase transition front is non-zero even for the symmetric potential, $\delta =0$. However, for zero kinetic modulus in the case of the symmetric potential the wave amplitude is zero.

As expected, for the constrained system we obtained solutions for the viscous CH equation (for non-zero kinetic modulus $\beta $), and for CH equation (for zero kinetic modulus $\beta $). It is easy to show that these solutions could be obtained directly for the system \eqref{2.6}-\eqref{2.7} after the formal limit $\varepsilon \to 0$.

The third exactly solvable case of the algebraic system \eqref{3.4}-\eqref{3.5}, \eqref{3.7}-\eqref{3.8} is the unconstraint case $\varepsilon \ne 0$ with zero viscosity $\eta $ (zero kinetic modulus $\beta $). Remarkably, to have the real value of the parameters it should be $\varepsilon >\varepsilon _{l} =\frac{4}{\left(\left|a_{1} \right|+a_{2} \right)^{2} } $, i.e., this case is strictly unconstraint even formally, different from the physical limitation for the static solution. Presuming $\varepsilon \gg \max \left[\varepsilon _{l} ,\, 1\right]$ we obtained two sets of the parameters, which define two solutions. The necessary condition for the non-negativity of limiting concentrations is the negativity of $\delta $, i.e. the definite asymmetry of the potential. The limiting values are equal for the both solutions. However, while for the solution with the higher velocity the concentration and the order parameter increase/decrease simultaneously, for the solution with lower velocity, if the order parameter $u$ increases with $z$, the concentration $c$ decreases, and vice versa. The concentration front is always `retreating'. Interesting, this behavior reminds the behavior of the exact coupled kink/anti-kink solution obtained by Rosenau [12] for the diffusion equation with $\frac{1}{u} $ diffusivity: the `retreating' kink and anti-kink merge and the solution disappears.

In the first and third exactly solvable cases our model allows to monitor the `degree of constraint' explicitly and observe its influence on the solution parameters. For example, Eqs.~\eqref{4.22}-\eqref{4.23} show that the velocity of the transition front increases with increasing constraint (increasing $\chi =\frac{1}{\varepsilon } $). Comparing the above considered cases, we see that in the unconstraint cases of Allen-Cahn type, i.e. for two independent variables, the system demonstrates broader set of possible regimes; in the constraint limit only CH process remains.

\setcounter{equation}{0}
\subsection*{Appendix 1} \label{s:A1}
\renewcommand{\theequation}{A1.\arabic{equation}}

%\textbf{\textit{Appendix 1}}

Using system \eqref{3.9}-\eqref{3.11} we calculate the values $u_{1s} ,\, u_{2s} $ of the order parameter, and the values of the concentration $c_{1s} ,\, c_{2s} $ at $\pm \infty $. Solving this system we get
\begin{eqnarray} \label{6.1} & & p_{s}=u_{1s}+u_{2s}=\frac{2}{3} \delta , \, \, \, \, \, \, q_{s}=u_{1s}u_{2s} =-\frac{2}{9} \delta ^{2} -s, \nonumber \\ & & \zeta _{s} =-\frac{1}{3} \delta \left(\frac{2}{9} \delta ^{2} +s\right) .  \end{eqnarray}
The values of the order parameter $u_{1s} ,\, u_{2s} $ at $\pm \infty $ for the stationary solution are the roots of the quadratic equation
\begin{equation} \label{6.2} u^{2} -p_{s} u+q_{s} =0 ,  \end{equation}
\begin{equation} \label{6.3} u_{1s,2s} =\frac{1}{3} \delta \mp \sqrt{\frac{1}{3} \delta ^{2} +s}  .  \end{equation}
Then the smaller value of the concentration is
\begin{equation} \label{6.4} c_{1s} =\phi _{0} \left(u_{1s} +\varepsilon \zeta _{s} \right) .  \end{equation}

While the order parameter could be negative, the concentration is definitely non-negative, so it should be $u_{1s} +\varepsilon \zeta _{s} \ge 0$, i.e.
\begin{equation} \label{6.5} \delta -\sqrt{3\delta ^{2} +9s} -\varepsilon \delta \left(\frac{2}{9} \delta ^{2} +s\right)\ge 0 .  \end{equation}
The latter inequality is definitely violated for $\delta \ge 0$ ($s$ is per definition positive). If $\delta <0$ Eq.~\eqref{6.5} becomes
\begin{equation} \label{6.6} \varepsilon \left|\delta \right|\left(\frac{2}{9} \delta ^{2} +s\right)\ge \sqrt{3\delta ^{2} +9s} +\left|\delta \right| .  \end{equation}

\setcounter{equation}{0}
\subsection*{Appendix 2} \label{s:A2}
\renewcommand{\theequation}{A2.\arabic{equation}}

%\textbf{\textit{Appendix 2}}

Here we check the non-negativity of the smaller concentration values $c_{11} ,\, \, c_{21} $. From equation \eqref{1.5} it follows
\begin{equation} \label{7.1} \frac{c}{\phi _{0} } =u+\varepsilon \mu  .\end{equation}
We remind that $\mu =\gamma u+\zeta $ and $\omega =\varepsilon \gamma +1$. Then \eqref{7.1} becomes
\begin{equation} \label{7.2} \frac{c}{\phi _{0} } =\omega u+\varepsilon \zeta . \end{equation}
For the solution to be physical the smallest concentration should be non-negative. If $\delta >0$, $\omega _{1} $ is positive and $\omega _{2} $ is negative, see \eqref{4.14}-\eqref{4.15}. So with $\omega _{1} $ we need to take in \eqref{7.2} the smallest `boundary' value $u_{11} $, and with $\omega _{2} $ we need to take the largest value $u_{22} $:
\begin{equation} \label{7.3} \frac{1}{\phi _{0} } c_{11} =\omega _{1} u_{11} +\varepsilon \zeta _{1} \ge 0, \end{equation}
\begin{equation} \label{7.4} \frac{1}{\phi _{0} } c_{21} =\omega _{2} u_{22} +\varepsilon \zeta _{2} \ge 0. \end{equation}

Now we consider these inequalities separately. The quantities $\zeta _{1} ,\, u_{11} $ are given by \eqref{4.16} and \eqref{4.19}; the inequality \eqref{7.3} is equivalent to
\begin{equation} \label{7.5} \omega _{1} u_{11} +\varepsilon \zeta _{1} =-\frac{2\varepsilon }{27} \delta \left(a_{2} +2\left|a_{1} \right|\right)\left[\left|a_{1} \right|+2a_{2} \right]\ge 0. \end{equation}
This inequality is definitely violated for positive $\delta $; so the solution with these values of parameters is unphysical. The quantities $\zeta _{2} ,\, u_{22} $ are given by \eqref{4.17} and \eqref{4.20}; inequality \eqref{7.4} is equivalent to
\begin{equation} \label{7.6} \omega _{2} u_{22} +\varepsilon \zeta _{2} =-\frac{2\varepsilon }{27} \delta \left(2a_{2} +\left|a_{1} \right|\right)\left(a_{2} +2\left|a_{1} \right|\right)\ge 0 .\end{equation}
This inequality coincides with \eqref{7.5} and is also violated for the positive $\delta $. On the other hand, for negative $\delta $ the parameter $\omega _{1} $ is negative, $\omega _{2} $ is positive, and the left hand sides of inequalities \eqref{7.5}-\eqref{7.6} are the largest and definitely positive values of non-dimensional concentrations, $c_{12} \sqrt{\rho } $ and $c_{22} \sqrt{\rho } $. Calculating the smallest values of concentrations for $\delta <0$, i.e.,
\begin{equation} \label{7.7} c_{11} \sqrt{\rho } =\omega _{1} u_{12} +\varepsilon \zeta _{1}, \end{equation}
\begin{equation} \label{7.8} c_{21} \sqrt{\rho } =\omega _{2} u_{21} +\varepsilon \zeta _{2},  \end{equation}
we get $c_{11} =c_{21} =0$.

{\small \topsep 0.6ex

}

\end{document}